\documentclass[12pt]{article}

\usepackage{authblk}
\usepackage{comment}
\usepackage{amssymb}
\usepackage{amsthm}
\usepackage{amscd}
\usepackage{graphicx}
\usepackage{titlesec}
\titleformat{\section}{\bfseries}{\thesection .}{0.5em}{}
\titleformat{\subsection}{\itshape}{\thesubsection .}{0.5em}{}
\usepackage[top=25.4mm, bottom=25.4mm, left=31.7mm, right=32.2mm]{geometry}
\usepackage[tbtags]{amsmath}
\usepackage{booktabs}
\allowdisplaybreaks  
\usepackage[,sort&compress]{natbib}

\usepackage{caption}
\usepackage{amsmath}
\usepackage{amsfonts}
\usepackage{graphicx,amsmath,amsfonts,amssymb,amscd,amsthm}
\usepackage{enumerate}
\newtheoremstyle{kai}
{3pt} {3pt} {} {} {\bfseries} {.} {.5em} {}
\def\EquationsBySection{\def\theequation
{\thesection.\arabic{equation}}%
\@addtoreset{equation}{section}}

\newcommand\old[1]{}
\newcommand{\pend}{\hfill \thicklines \framebox(6.6,6.6)[l]{}}
\renewenvironment{proof}{\noindent {\it  Proof.} \rm}{\pend}
\newtheorem{theorem}{Theorem}[section]

\newtheorem{corollary}{Corollary}[section]
\newtheorem{proposition}{Proposition}[section]
\newtheorem{remark}{Remark}[section]
\newtheorem{definition}{Definition}[section]

\newtheorem{algorithm}{Algorithm}[section]
\EquationsBySection \makeatother
\usepackage{indentfirst}

\begin{document}
\pagestyle{plain}
\title
{\bf On the $BMAP_1, BMAP_2/PH/g, c$ retrial queueing system}
\author{Jinbiao Wu $^a$\thanks{E-mail address: wujinbiao@ymail.com}, {Yi
Peng$^b$}, {Zaiming Liu$^a$}
\\
\vspace{0.2cm}\small\it $^a$School of Mathematics and Statistics, \\
 \small\it Central South University,
Changsha 410083, Hunan, P.R.  China\\
 \small\it  $^b$Junior Education Department,\\ \small\it Changsha Normal University, Changsha
410100,  Hunan, P.R.  China}

\date{}
\maketitle

\begin{abstract}
In this paper, we analyze a retrial queueing system with Batch
Markovian Arrival Processes and two types of customers. The rate of
individual repeated attempts from the orbit is modulated according
to a Markov Modulated Poisson Process. Using the theory of
multi-dimensional asymptotically quasi-Toeplitz Markov chain, we
obtain the stability condition and the algorithm for calculating the
stationary state distribution of the system. Main performance
measures are presented. Furthermore, we investigate some
optimization problems. The algorithm for determining the optimal
number of guard servers and total servers is elaborated. Finally,
this queueing system is applied to the cellular wireless network.
Numerical results to illustrate the optimization problems and the
impact of retrial on performance measures are provided. We find that
the performance measures are mainly affected by the two types of
customers' arrivals and service patterns, but the retrial rate plays
a less crucial role.
\end{abstract}

\noindent{\it Keywords:} \small  Queueing; Markov Modulated Poisson
Process; Asymptotically quasi-Toeplitz Markov chains; Cellular
wireless networks; Performance evaluation

\section{Introduction}\label{sec1}
Most queueing systems assume that the arrival process is a
stationary Poisson process. But such a process does not characterize
the typical features of traffic in modern telecommunication networks
such as correlation and burstiness. The batch Markovian arrival
process (BMAP) is a useful and appropriate model for describing such
features. The BMAP is a generalization of many well-known processes
including the Markovian arrival process (MAP), the Markov-modulated
Poisson process (MMPP), the PH-renewal process and the Poisson
process. Furthermore, the BMAP preserves the tractable Markovian
structure. The origins of the BMAP can be traced to the development
of the versatile Markovian point process by  \cite{Neuts1979}. The
currently used denotations were introduced in  \cite{Neuts1989}. It
is widely employed in research for queueing systems, inventory
systems, reliability engineering, manufacturing systems, computer
communication systems, insurance problems and so on.
\cite{Chakravarthy1999} gave a detailed survey of the results about
the queues which the input flow is BMAP. \cite{Heyman} presented an
application of the BMAP for modeling the information flows in the
modern telecommunication networks.

Retrial queues are very good mathematical models for cellular
wireless networks, telephone switching systems, local area networks
under the protocols of random multiple access, computer systems for
competing to gain service from a central processor unit, etc. Over
recent years it has been a rapid growth in the literature on retrial
queues. The reader can refer to the survey papers by \cite{Falin},
 \cite{Artalejo1,Artalejo2}, \cite{Yang} for a review of
main results and methods. Most of the papers and books are focused
on single server retrial queues. Multi-server retrial queues are
more complicated for research and even the model M/M/c is studied
until now only algorithmically via different truncation schemes. The
multi-server retrial queues with the BMAP are still not well
investigated in literature. In \cite{Breuer2002}, the classical
BMAP/PH/N retrial queue was analyzed by means of discrete-time
multi-dimensional asymptotically quasi-Toeplitz Markov chain.
\cite{He2000} gave the stable condition for the BMAP/PH/S/S+K
retrial queue with PH-retrial times. Recently, \cite{Klimenok2007}
analyzed the BMAP/PH/N retrial queue with impatient customers.
Later, \cite{Kim2008} studied the BMAP/PH/N retrial queue with
Markovian flow of breakdowns. Subsequently, \cite{Kim2010}
considered the BMAP/PH/N retrial queueing system operating in
Markovian random environment.  \cite{Dudin2012} investigated the
BMAP/PH/N retrial queue with Markov modulated retrials.
Particularly, for the retrial queues with two types of customers,
the single server retrial queues have received considerable
attention in the literature. For example, using the supplementary
variable method,  \cite{Choi1990} dealt with an $M_1,M_2/G/1$
retrial queue with two types of customers and non-preemptive resume
priority. After that, Choi and Park's model was extended to an
$M_1,M_2/G_1,G_2/1$ retrial queue with priority customers by
\cite{Falin1993}. Later,  \cite{Choi1995} investigated an $M/G/1$
retrial queue with two types of customers and finite capacity.
Furthermore, using the matrix analytic method, \cite{Martin} studied
an $M_1,M_2/G_1,G_2/1/1$ retrial queue with impatient customers. In
addition,  \cite{Choi1999} presented a survey of single server
retrial queue with two types of customers and priority. Recently,
\cite{Wu2013} introduced negative customers into the retrial queues
with two types of customers. However, to the best of our knowledge,
little work has been done on multi-server retrial queues with two
types of customers, and no work on BMAP/PH/N retrial queue with two
types of customers is found in the queueing literature.

Our model can be successfully studied by means of asymptotically
quasi-Toeplitz Markov chains. Here, we briefly introduce the
continuous-time multi-dimensional asymptotically quasi-Toeplitz
Markov chains. For a comprehensive and excellent overview of
multi-dimensional asymptotically quasi-Toeplitz Markov chains with
discrete and continuous time, readers may refer to the paper
\cite{Klimenok2006}. The next definition and propositions follow
 \cite{Klimenok2006}. Let $\xi_t=\{i_t,
\mathbf{a}_t\}$,  $t\geq0$, be a regular irreducible continuous time
Markov chain with the state space $\Omega=\{(i, \mathbf{a}),
\mathbf{a}\in L_i, i=0,1, \cdots, i^0; (i, \mathbf{a}),
\mathbf{a}\in L, i>i^0\}$. Enumerate the states of the chain
$\{\xi_t$, $t\geq0\}$ as follows: the states $(i, \mathbf{a})$ are
numerated in ascending order of the component $i$ and for fixed $i$,
the states $(i, \mathbf{a})$ are arranged in lexicographic order. We
denote the block matrix $Q=(Q_{i,j})$ as the generator of the chain,
where the block $Q_{i,j}$ is formed by intensities
$q_{(i,\mathbf{a});(j,\mathbf{b})}$ of the chain, transition from
the state $(i,\mathbf{a})$ to the state $(j,\mathbf{b})$.
\begin{definition}
A regular irreducible continuous time Markov chain $\{\xi_t,
t\geq0\}$ is called asymptotically quasi-Toeplitz Markov chain if
\begin{enumerate}[(i)]
\item $Q_{i,j}=0$ for $j<i-1$, $i>0$.
\item There exit matrices $Y_n, n\geq0$, such that
\begin{align*}
Y_n&=\lim_{i\to\infty}\Lambda_i^{-1}Q_{i, i+n-1}, \ \ \ n=0, 2,3,
\cdots,\\
Y_1&=\lim_{i\to\infty}\Lambda_i^{-1}Q_{i, i}+I,
\end{align*}
and the matrix $\sum_{n=0}^\infty Y_n$ is stochastic.
\item The jump chain $\{\xi_n, n\geq1\}$ of the process $\{\xi_t,
t\geq0\}$ is non-periodic,
\end{enumerate}
where the diagonal matrix $\Lambda_i$ defines the diagonal entries
of the block $-Q_{ii}$, $i\geq0$.
\end{definition}

Let $Y(z)=\sum_{k=0}^\infty Y_kz^k$ be the matrix generating
function characterizing the limiting Markov chain. The following two
propositions present the ergodicity conditions for the Markov chain
$\{\xi_t, t\geq0\}$ in terms of the generating function $Y(z)$ based
on whether the matrix $Y(1)$ is an irreducible matrix or not.
\begin{proposition}
If $Y(1)$ is an irreducible matrix and assume that
\begin{enumerate}[(i)]
\item the series $\sum_{n=1}^\infty n Q_{i,i+n-1}\mathbf{e}$
converges for $i=\overline{0, i^0}$,
\item the series $\sum_{n=1}^\infty n Q_{i,i+n-1}\mathbf{e}$
converges for $i>i^0$, and there exists a positive integer $i^*>i^0$
such that this series converges uniformly in the region $i\geq i^*$,
\item the sequence $\Lambda_i^{-1}$, $i\geq0$, has an upper bound,
\end{enumerate}
then the sufficient condition for the Markov chain $\{\xi_t,
t\geq0\}$ ergodicity is the fulfillment of the inequality
\begin{align*}
\left[\det(zI-Y(z))\right]'|_{z=1}>0.
\end{align*}
\end{proposition}

If $Y(1)$ is a reducible matrix, then it can be presented in a
normal form which has the following structure:
\begin{align*} Y(z)=\left(
\begin{array}{ccccc}
Y^{(1)}(z)& & & & \\
 &Y^{(2)}(z)& & & \\
 & &\ddots& & \\
 & & &Y^{(l)}(z)& \\
 Y^{(l+1,1)}(z)& Y^{(l+1,2)}(z)& \cdots &Y^{(l+1,l)}(z)&Y^{(l+1)}(z)
\end{array}
\right),
\end{align*}
where $Y^{(i)}(z)$, $i=\overline{1, l}$, are irreducible square
matrices.
\begin{proposition}
If $Y(1)$ is a reducible matrix of the above form and the conditions
of (i)-(iii) of Proposition 1.1 are satisfied, then the sufficient
condition for the Markov chain $\{\xi_t, t\geq0\}$ ergodicity is the
fulfillment of the inequality
\begin{align*}
\left[\det(zI-Y^{(i)}(z))\right]'|_{z=1}>0, \ \ \ i=\overline{1, l}.
\end{align*}
\end{proposition}

The rest of the paper is organized as follows. The model description
is given in Section 2. The Markov chain associated with the system
is analyzed in Section \ref{sec3} where we formulate the model as a
multi-dimensional asymptotically quasi-Toeplitz Markov chain and
obtain the ergodicity condition of the Markov chain. The algorithms
for calculating the stationary state probabilities are elaborated in
Section 4. A set of performance measures is demonstrated in Section
5. Two optimization problems and an algorithm for calculating the
optimal values are provided in Section 6. An application for the
model under discussion and some numerical examples are shown in
Section 7. Finally, summary of the results is presented in
Conclusion.

 \section{Model Description}\label{sec2}

 We consider a multi-server retrial queueing system. Two types of customers (primary customers
and priority customers) arrive to the system according to two
independent BMAPs (denoted by $BMAP_1$ and $BMAP_2$). The notion of
the BMAP and its detailed description were given by
\cite{Neuts1989}. We denote the directing process of the $BMAP_1$ by
$\{\nu_t, t\geq0\}$ with the state space $\{1,\cdots,W\}$ and the
directing process of the $BMAP_2$ by $\{\gamma_t, t\geq0\}$ with the
state space $\{1,\cdots,V\}$. The two directing processes behave as
two independent irreducible continuous-time Markov chains. Denote
the matrix sequences $(D_n: n\in \mathbb{N}_0)$ and $(E_n: n\in
\mathbb{N}_0)$ as the sequences of characterizing matrices for the
$BMAP_1$ and $BMAP_2$, respectively. Introduce the matrix generating
functions $D(z)=\sum_{n=0}^\infty D_n$ and $E(z)=\sum_{n=0}^\infty
E_n$, $|z|<1$. We assume the matrix functions $D(z)$ and $E(z)$
satisfy all assumptions of \cite{Lucantoni}. Therefore, the matrix
$D(1)$ is the generator of the process $\{\nu_t, t\geq0\}$ and the
matrix $E(1)$ is the generator of the process $\{\gamma_t,
t\geq0\}$. Suppose that the two arrival processes $BMAP_1$ and
$BMAP_2$ start with the initial phase probability distributions
$\alpha_1$ and $\alpha_2$, respectively. Let $\theta_1$ and
$\theta_2$ be the stationary probability vectors of $D(1)$ and
$E(1)$, respectively. Then $\lambda_1=\theta_1 D'(1)\mathbf{e}$ is
the stationary arrival rate of the primary customers and
$\lambda_2=\theta_2 E'(1)\mathbf{e}$ is the stationary arrival rate
of the priority customers. Moreover, we define
$\lambda_{b1}=\theta_1(-D_0)\mathbf{e}$ and
$\lambda_{b2}=\theta_2(-E_0)\mathbf{e}$ as the intensities of group
primary arrivals and group priority arrivals, respectively. Here
$\mathbf{e}$ is a column vector of appropriate size with all
elements equal to 1.

The service facility consists of $c$ ($c\geq2$) identical servers
which are identical and operate independently of each other. All
customers have the same independent, identical PH distributed
service times $S$ which is governed by the irreducible continuous
time Markov chain $\{m_t, t\geq0\}$ with the state space
$\{1,\cdots,M\}$. The transitions of the Markov chain $\{m_t,
t\geq0\}$, which do not lead to service completion, are defined by
the irreducible matrix $S$ of size $M\times M$. The transitions of
the Markov chain $\{m_t, t\geq0\}$, which lead to service
completion, are defined by the vector $S_0=-S\mathbf{e}$. At the
service beginning epoch, the state of the process $\{m_t, t\geq0\}$,
is determined by the probabilistic row vector $\varsigma$ of size
$1\times M$. The service time distribution function has the form
$B(x)=1-\varsigma e^{Sx}\mathbf{e}$. The mean rate of service is
$\mu=[-\varsigma S^{-1}\mathbf{e}]^{-1}$. A more detailed
information about the PH type distribution can be seen in the books
 \cite{Neuts1981},  \cite{Latouche} and
 \cite{He2014}.

When a batch of priority customers arrives at the system and there
are several servers being idle, the priority customers occupy the
corresponding number of the servers. If the number of the idle
servers is insufficient (or all the servers are busy) the rest of
the batch (or all the batch) enters the retrial group (called
orbit). When a batch of primary customers arrives at the system and
meets no more than $g-1$ ($1\leq g\leq c-1$) busy servers, the
primary customers occupy the corresponding number of the servers (If
the number of the idle servers is insufficient, the rest of the
batch goes to the orbit); otherwise, the whole batch of primary
customers is blocked and goes to the orbit. These customers in orbit
are said to be repeated customers and the orbit capacity is assumed
to be unlimited. For all customers in the orbit, the individual
repeated attempts are governed by a common Markov Modulated Poisson
Process (MMPP). The directing process of the Markov Modulated
Poisson Process is denoted by $\{r_t, t\geq0\}$ which is an
irreducible continuous time Markov chain with the state space
$\{1,\cdots,R\}$. Transition intensities of the Markov chain $\{r_t,
t\geq0\}$, which are not related with customers retrials, are
defined by the matrix $T_0$, and transition intensities, which are
accompanied by retrial, are described by the  diagonal matrix
$T_1=\mbox{diag}\{\sigma_1, \cdots, \sigma_R\}$. The matrices
$T=T_0+T_1$ is the infinitesimal generators of the processes $\{r_t,
t\geq0\}$. That is, when the process $r_t$ stays in the state $r$,
$r=1, \cdots, R$, each customer from the orbit makes attempts to
seek the service  at exponentially distributed time intervals with
mean $1/\sigma_r$. Hence, the total rate of retrials is equal to
$n\sigma_r, \sigma_r>0$, when the process $r_t$ is in the state $r$
and the orbit size (the number of calls on the orbit) is equal to
$n$, $n>0$. Let $\theta_0$ be the stationary probability vector of
$T$. Then $\sigma=\theta_0 T_1\mathbf{e}$ is the stationary retrial
rate of the repeated calls in the orbit. On retrial, an orbiting
customer obtains service immediately from one of the idle servers if
the number of busy servers is no more than $g-1$; otherwise it will
return to the orbit for later retrial. We understand that there are
$c-g$ servers (called guard servers) opening only for the priority
customers.

\indent For the use in the sequel, let us introduce the following
notations:\\
$\star$  $\mathbf{e}$ is a column vector of ones of suitable size.
When
needed we will identify the dimension of this vector with a subindex;\\
$\star$ $I$ is an identity matrix of appropriate size. When needed
we will identify the dimension of this matrix with a subindex;\\
$\star$ $O_L$ is a zero matrix of size $L\times L$;\\
$\star$ $\mbox{\emph{diag}}\{a_r, r=\overline{1,L}\}$ is a diagonal
matrix with diagonal entries $a_r$, $r=\overline{1,L}$;\\
$\star$  $\otimes$ and $\oplus$ are the symbols of Kronecker product
and sum of matrices, see, e.g. \cite{Graham};\\
$\star$  $\Omega^{\otimes
l}=\underbrace{\Omega\otimes\cdots\otimes\Omega}\limits_l$,
$l\geq1$,
 $\Omega^{\otimes0}=1, \Omega^{\oplus0}=0$;\\
$\star$  $\Omega^{\oplus
l}=\sum\limits_{m=0}^{l-1}I_{n^m}\otimes\Omega\otimes
I_{n^{l-m-1}}$, $l\geq1$, for the matrix $\Omega$ having $n$ rows.

Our work objectives are to derive the stability condition, to
calculate the stationary state distribution, to derive the main
performance measures of the system and optimize the values of $g$
and $c$.

\section{Analysis}\label{sec3}
In this section, we first adopt a layered approach for the
description of the state space to construct the infinitesimal
generator $\mathbf{Q}$ of the Markov chain which can represent the
model. We then derive the ergodicity condition of the Markov chain.

Let\\
$\ast$  $o_t$ be the number of repeated customers in the orbit, at
the epoch $t$, $t\geq0$,
$o_t\geq0$;\\
$\ast$  $b_t$ be the number of busy servers, at the epoch $t$,
$t\geq0$,
$b_t=\overline{0,c}$;\\
 $\ast$ $r_t$ be the state of the directing process of the MMPP
of repeated customers,  at the epoch $t$,
$t\geq0$, $r_t=\overline{1,R}$;\\
 $\ast$ $\nu_t$ be
the state of the directing process of the BMAP of primary customers,
at the epoch $t$, $t\geq0$, $\nu_t=\overline{1,W}$;\\
$\ast$ $\gamma_t$ be the state of the directing process of the BMAP
of priority customers, at the epoch $t$, $t\geq0$, $\gamma_t=\overline{1,V}$;\\
$\ast$ $m_t^{(j)}$ be the state of the directing process of the
service on the $j$th busy server, at the epoch $t$, $t\geq0$,
$m_t^{(j)}=\overline{1,M}$, $j=\overline{1,b_t}$. Here, we assume
that the busy servers are numerated in order of their occupying,
i.e., the server, which begins the service, is appointed the maximal
number among all busy servers; when some server finishes the service
and is released, the servers are correspondingly enumerated.

\indent Then, the process of the model can be described by the
regular irreducible continuous time multi-dimensional Markov chain:
$$\xi=\{\xi_t, t\geq0\}=\left\{o_t,b_t, r_t, \nu_t, \gamma_t, m_t^{(1)},\cdots,m_t^{(b_t)}, t\geq
0\right\},$$ with the state space
\begin{eqnarray*}
\Omega&=&\left\{\left(o, b, r, \nu,\gamma,m^{(1)},\cdots,
m^{(b)}\right);\right.\\  & & \left. o\geq0, b=\overline{0,c},
r=\overline{1,R}, \nu=\overline{1,W}, \gamma=\overline{1,V},
m^{(i)}=\overline{1,M}, i=\overline{1,b}\right\},
\end{eqnarray*}
and the dimension of the state space of the Markov chain $\{\xi_t,
t\geq0\}$ is equal to $K=RWV\frac{1-M^{c+1}}{1-M}$. We note that the
dimension $K$ can be very large. E.g., if the array $(c,R,W,V,M)$ is
equal to $(10,2,2,2,2)$, then $K=16376$.

 \indent  We suppose that the states of the Markov chain $\{\xi_t,
t\geq0\}$ are enumerated in lexicographic order. By means of
considering the probabilities of the Markov chain transitions during
an infinitesimal time interval, the infinitesimal generator matrix
$\mathbf{Q}$ of the Markov chain $\{\xi_t, t\geq0\}$ has the
following structure:
\begin{align} \mathbf{Q}=\left(
\begin{array}{ccccc}
\mathbf{Q}_{00}&\mathbf{Q}_{01}& \mathbf{Q}_{02}& \mathbf{Q}_{03}& \cdots\\
\mathbf{Q}_{10}&\mathbf{Q}_{11}&\mathbf{Q}_{12}&\mathbf{Q}_{13} & \cdots\\
 &\mathbf{Q}_{21}&\mathbf{Q}_{22}&\mathbf{Q}_{23}& \cdots\\
 & &\mathbf{Q}_{32}&\mathbf{Q}_{33}&\cdots\\
 & & &\ddots&\ddots
\end{array}
\right),
\end{align}

where the blocks $\mathbf{Q}_{ij}$ of size $K\times K$ are defined
by:
\begin{align*}
(\mathbf{Q}_{i,i})_{l,l'}=\left\{
\begin{array}{ll}
O, & l'<l-1, l=\overline{2, c},\\
I_{RWV}\otimes S_0^{\oplus l},
& l'=l-1, l=\overline{1, c},\\
T_0\oplus D_0\oplus E_0\oplus S^{\oplus
l}-iT_1 \otimes I_{WVM^l},&l'=l, l=\overline{0, g-1},\\
T_0\oplus D_0\oplus E_0\oplus S^{\oplus l},& l'=l, l=\overline{g, c},\\
I_R\otimes(D_r\oplus E_r)\otimes I_{M^l}\otimes \varsigma^{\otimes r}, & l'=l+r,
r=\overline{1, g-l}, l=\overline{0, g-1},\\
I_{RW}\otimes E_r\otimes I_{M^l}\otimes\varsigma^{\otimes r}, &
l'=l+r,  r=\overline{g+1-l, c-l}, l=\overline{0, g-1},\\
 I_{RW}\otimes
E_r\otimes I_{M^l}\otimes\varsigma^{\otimes r}, & l'=l+r,
r=\overline{1, c-l}, l=\overline{g, c-1}.
\end{array}
\right. \ \ \ i\geq0,
\end{align*}

 {\scriptsize
\begin{align*}
       \mathbf{Q}_{i,i-1}=i
         \bordermatrix{
         & \scriptstyle0 & \scriptstyle 1 & \scriptstyle2 & \cdots & \scriptstyle g &  \cdots &\scriptstyle c \cr
         \scriptstyle0 &  &  T_1\otimes I_{WV}\otimes\varsigma& & & & & \cr
         \scriptstyle 1 &  &   & T_1\otimes I_{WVM}\otimes\varsigma& & & & \cr
         \vdots & & & & \ddots& & & \cr
         \scriptstyle g-1 & & & & &T_1\otimes I_{WVM^{g-1}}\otimes\varsigma & & \cr
         \scriptstyle g & & & & & & & \cr
         \vdots &   & & & & & & \cr
         \scriptstyle c &  & & & & & & },\ \ \ i\geq1,
\end{align*}
\begin{align*}
       \mathbf{Q}_{i,i+k}=
         \bordermatrix{
         & \scriptstyle0 & \scriptstyle 1 & \cdots & \cdots & \scriptstyle g &  \cdots &\scriptstyle c-1&\scriptstyle c \cr
         \scriptstyle0 &  & & & & I_R\otimes D_{k+g}\otimes I_{V}\otimes \varsigma^{\otimes g}& & & I_{RW}\otimes E_{k+c}\otimes \varsigma^{\otimes c}\cr
         \scriptstyle 1 &  & & & & I_R\otimes D_{k+g-1}\otimes I_{VM}\otimes \varsigma^{\otimes (g-1)}& & & I_{RW}\otimes E_{k+c-1}\otimes  I_{M}\otimes \varsigma^{\otimes (c-1)}\cr
         \vdots & & & & & \vdots& & & \vdots\cr
         \scriptstyle g-1 & & & & &  I_R\otimes D_{k+1}\otimes I_{VM^{g-1}}\otimes \varsigma& & & \vdots\cr
         \scriptstyle g & & & & &I_R\otimes D_{k}\otimes I_{VM^{g}} & & & \vdots\cr
         \vdots &   & & & & &\ddots & & \vdots\cr
         \scriptstyle c-1 &   & & & & & &I_R\otimes D_k\otimes I_{VM^{c-1}}  & I_{RW}\otimes E_{k+1}\otimes  I_{M^{c-1}}\otimes \varsigma\cr
         \scriptstyle c &  & & & &  & & &I_R\otimes D_k\otimes I_{VM^{c}}+I_{RW}\otimes E_k\otimes  I_{M^{c}}},
\end{align*}
} $i\geq0$, $k\geq1$.

\begin{remark}
When $E_n=0$, $n=0,1, \cdots$ and $c=g$, the present system reduces
to a retrial BMAP/PH/c queueing system with Markov modulated
retrials. It may be noted that the results above after putting
$E_n=0$, $n=0,1, \cdots$ and $c=g$ agree with the result presented
in  \cite{Dudin2012}.
\end{remark}

In the following sections, we analyze the queueing system by using
the theory of multi-dimensional asymptotically quasi-Toeplitz Markov
chains.

 The diagonal blocks of the generator $\mathbf{Q}$ are given
by
\begin{align*}
(Q_{ii})_{l,l}=\left\{
\begin{array}{ll}
T_0\oplus D_0\oplus E_0\oplus S^{\oplus
l}-iT \otimes I_{WVM^l},&l=\overline{0, g-1},\\
T_0\oplus D_0\oplus E_0\oplus S^{\oplus l},&  l=\overline{g, c},
\end{array}
\right. , \ \ \ i\geq0.
\end{align*}

Denote the diagonal entries of the matrices $D_0$, $E_0$ and $S$ as
$\{-\lambda^{(1)}_w, w=\overline{1, W}\}$,  $\{-\lambda^{(2)}_v,
v=\overline{1, V}\}$ and $\{-s_m, m=\overline{1, M}\}$ respectively.
Then the diagonal matrix $\Lambda_i$ in Definition 1.1 is given by:
\begin{align}
\Lambda_i=\Delta+i\tilde{T}_1\hat{I},
\end{align}
where \begin{align*}
\Delta=&\mbox{\emph{diag}}\{\mbox{\emph{diag}}\{\sigma_r,
r=\overline{1,
R}\}\oplus \mbox{\emph{diag}}\{\lambda^{(1)}_w, w=\overline{1, W}\} \\
&\oplus \mbox{\emph{diag}}\{\lambda^{(2)}_v, v=\overline{1,
V}\}\oplus [\mbox{\emph{diag}}\{s_m, m=\overline{1, M}\}]^{\oplus
l}, l=\overline{0, c}\},
\end{align*}
\begin{align*}
\tilde{T}_1=\mbox{\emph{diag}}\{T_1\otimes I_{WVM^l}, l=\overline{0,
c}\},
\end{align*}
\begin{align*}
\hat{I}=\left(
\begin{array}{cc}
I_{RWV\sum_{n=0}^{g-1}M^n}&  \\
 & O_{RWV\sum_{n=g}^cM^n}
\end{array}
\right).
\end{align*}

Further, introduce the matrices $\bar{I}$, $\Upsilon$ and
$\Gamma_k$, $k\geq0$, of size $K\times K$:
\begin{align*}
\bar{I}=I-\hat{I}=\left(
\begin{array}{cc}
O_{RWV\sum_{n=0}^{g-1}M^n}&  \\
 & I_{RWV\sum_{n=g}^cM^n}
\end{array}
\right),
\end{align*}
\begin{align*}
       \Upsilon=
         \bordermatrix{
         & \scriptstyle0 & \scriptstyle 1 & \scriptstyle2 & \cdots & \scriptstyle g &  \cdots &\scriptstyle c \cr
         \scriptstyle0 & O_{RWV} &  I_{RWV}\otimes\varsigma& & & & & \cr
         \scriptstyle 1 &  & O_{RWVM}  & I_{RWVM}\otimes\varsigma& & & & \cr
         \vdots & & & \ddots& \ddots& & & \cr
         \scriptstyle g-1 & & & & \ddots&I_{RWVM^{g-1}}\otimes\varsigma & & \cr
         \scriptstyle g   & & & &       & O_{RWVM^g}& & \cr
         \vdots &   & & & & & \ddots& \cr
         \scriptstyle c &  & & & & & & O_{RWVM^c}},
\end{align*}

\begin{align*}
(\Gamma_0)_{l,l'}=\left\{
\begin{array}{ll}
O, & l'<l-1, l=\overline{2, c},\\
I_{RWV}\otimes S_0^{\oplus l},
& l'=l-1, l=\overline{1, c},\\
T_0\oplus D_0\oplus E_0\oplus S^{\oplus
l},&l'=l, l=\overline{0, c},\\
I_R\otimes(D_r\oplus E_r)\otimes I_{M^l}\otimes \varsigma^{\otimes r}, & l'=l+r, r=\overline{1, g-l},
l=\overline{0, g-1},\\
I_{RW}\otimes E_r\otimes I_{M^l}\otimes\varsigma^{\otimes r}, &
l'=l+r,  r=\overline{g+1-l, c-l}, l=\overline{0, g-1},\\
 I_{RW}\otimes
E_r\otimes I_{M^l}\otimes\varsigma^{\otimes r}, & l'=l+r,
r=\overline{1, c-l}, l=\overline{g, c-1}.
\end{array}
\right.
\end{align*}
$$\Gamma_k=Q_{i,i+k},  \ \ \ k\geq1.$$
We obtain:
\begin{align}\label{33}
\lim_{i\to\infty}\Lambda_i^{-1}=\Delta^{-1}\bar{I},
\end{align}
\begin{align}\label{34}
Y_0=\lim_{i\to\infty}\Lambda_i^{-1}Q_{i,i-1}= \Upsilon,
\end{align}
\begin{align}\label{35}
Y_1=\lim_{i\to\infty}\Lambda_i^{-1}Q_{ii}+I=\Delta^{-1}\bar{I}\Gamma_0+\bar{I},
\end{align}
\begin{align}\label{36}
Y_k=\lim_{i\to\infty}\Lambda_i^{-1}Q_{i,i+k-1}=\Delta^{-1}\bar{I}\Gamma_{k-1},
\ \ \ k\geq2,
\end{align}
and
\begin{align}\label{37}
\sum_{k=0}^\infty Y_k \mathbf{e}=\mathbf{e}.
\end{align}

 It can be easily verified that the Markov chain $\{\xi_t,
t\geq0\}$ satisfies Definition 1.1. Hence, it belongs to the class
of asymptotically quasi-Toeplitz Markov chain. From
\eqref{34}-\eqref{36}, we get
\begin{align}\label{37}
Y(z)=\Upsilon+\bar{I}z+\Delta^{-1}\bar{I}\Gamma(z)z, \ \ \ |z|<1,
\end{align}
where
\begin{align*}
(\Gamma(z))_{l,l'}&=\left(\sum_{k=0}^\infty\Gamma_kz^k\right)_{l,l'}\\
&=\left\{
\begin{array}{ll}
O, & l'<l-1, l=\overline{2, c},\\
I_{RWV}\otimes S_0^{\oplus l},
& l'=l-1, l=\overline{1, c},\\
T_0\oplus D_0\oplus E_0\oplus S^{\oplus
l},&l'=l, l=\overline{0, g-1},\\
\Xi_{g-l}(z), &l'=g, l=\overline{0, g-1},\\
\Theta_l(z), &l'=l, l=\overline{g, c},\\
I_R\otimes(D_r\oplus E_r)\otimes I_{M^l}\otimes \varsigma^{\otimes
r}, & l'=l+r, r=\overline{1, g-1-l},
l=\overline{0, g-2},\\
I_{RW}\otimes E_r\otimes I_{M^l}\otimes\varsigma^{\otimes r}, &
l'=l+r,  r=\overline{g+1-l, c-1-l}, l=\overline{0, g-1},\\
 I_{RW}\otimes
E_r\otimes I_{M^l}\otimes\varsigma^{\otimes r}, & l'=l+r,
r=\overline{1, c-1-l}, l=\overline{g, c-2},\\
\Psi_{c-l}(z), & l'=c, l=\overline{0, c-1},
\end{array}
\right.
\end{align*}
here
\begin{align*}
\Xi_n(z)=I_R\otimes\left[I_W\otimes
E_n+z^{-n}\left(D(z)-\sum_{k=0}^{n-1}D_kz^k\right)\otimes
I_{VM^{g-n}}\right]\otimes \varsigma^{\otimes n}, \ \ \
n=\overline{1, g},
\end{align*}
\begin{align*}
\Theta_n(z)=T_0\oplus D(z)\oplus E_0\oplus S^{\oplus n}, \ \ \
n=\overline{g, c-1},
\end{align*}
\begin{align*}
\Theta_c(z)=T_0\oplus D(z)\oplus E(z)\oplus S^{\oplus c},
\end{align*}
\begin{align*}
\Psi_n(z)=I_{RW}\otimes
z^{-n}\left(E(z)-\sum_{k=0}^{n-1}E_kz^k\right)\otimes
I_{M^{c-n}}\otimes \varsigma^{\otimes n}, \ \ \ n=\overline{1, c}.
\end{align*}

Substituting the expressions of the matrices $\Upsilon$, $\bar{I}$,
$\Delta$ and $\Gamma(z)$ in to \eqref{37}, we see that the matrix
function $\Gamma(z)$ is of the form
\begin{align}\label{38}
Y(z)=\left(
\begin{array}{cc}
Y_{11}&Y_{12}  \\
 & Y_{22}(z)
\end{array}
\right),
\end{align}
where
\begin{align*}
       Y_{11}=
         \bordermatrix{
         & \scriptstyle0 & \scriptstyle 1 & \scriptstyle2 & \cdots & \scriptstyle g-2  \cr
         \scriptstyle0 &O_{RWV}  &  I_{RWV}\otimes\varsigma& & &  \cr
         \scriptstyle 1 &  & \ddots  & I_{RWVM}\otimes\varsigma& &  \cr
         \vdots & & & \ddots& \ddots&  \cr
         \scriptstyle g-3 & & & &\ddots &I_{RWVM^{g-3}}\otimes\varsigma  \cr
          \scriptstyle g-2 & & & & & O_{RWVM^{g-2}}}_{RWV\sum_{i=0}^{g-2}M^i},
\end{align*}
\begin{align*}
       Y_{12}=
         \bordermatrix{
         &  &  &  &  &   \cr
         \scriptstyle0 &  &  & & &  \cr
         \scriptstyle 1 &  &   & & &  \cr
         \vdots & & & & &  \cr
          \scriptstyle g-2 & I_{RWVM^{g-2}}\otimes\varsigma& & & & }_{RWV\sum_{i=0}^{g-2}M^i\times
          RWV\sum_{i=g-1}^{c}M^i},
\end{align*}
\begin{align*}
       Y_{22}(z)=\left(
\begin{array}{cc}
O_{RWVM^{g-1}}& Y_{22}^{12}\\
 Y_{22}^{21}(z)& Y_{22}^{22}(z)
\end{array}
\right),
\end{align*}
\begin{align*}
       Y_{22}^{12}=\left(I_{RWVM^{g-1}}\otimes\varsigma,
       O_{RWVM^{g-1}\times RWV\sum_{n=(g+1)}^cM^n}\right),
\end{align*}
\begin{align*}
       Y_{22}^{21}(z)=\left(z\Delta_{gg}^{-1}\cdot I_{RWV}\otimes H_0^{\oplus g},
       O_{RWVM^{g-1}\times RWV\sum_{n=(g+1)}^cM^n}\right)^T,
\end{align*}
 {\scriptsize
\begin{align*}
Y_{22}^{22}(z)=z\left(
\begin{array}{ccccc}
\Delta_{gg}^{-1}\Theta_g(z)+I &\Delta_{gg}^{-1}\Gamma_{g,g+1} & \cdots & \Delta_{gg}^{-1}\Gamma_{g,c-1}& \Delta_{gg}^{-1}\Psi_{c-g}(z)\\
\Delta_{g+1,g+1}^{-1}\cdot I_{RWV}\otimes H_0^{\oplus g+1}&\Delta_{g+1,g+1}^{-1}\Theta_{g+1}(z)+I&\cdots&\Delta_{g+1,g+1}^{-1}\Gamma_{g+1,c-1} & \Delta_{g+1,g+1}^{-1}\Psi_{c-g-1}(z)\\
 &\Delta_{g+2,g+2}^{-1}\cdot I_{RWV}\otimes H_0^{\oplus g+2}&\cdots&\Delta_{g+2,g+2}^{-1}\Gamma_{g+2,c-1}& \Delta_{g+2,g+2}^{-1}\Psi_{c-g-2}(z)\\
 & &\ddots&\vdots&\vdots\\
  & &
  &\Delta_{c-1,c-1}^{-1}\Theta_{c-1}(z)+I&\Delta_{c-1,c-1}^{-1}\Psi_{1}(z)\\
 & & &\Delta_{c,c}^{-1}\cdot I_{RWV}\otimes H_0^{\oplus c}&\Delta_{c,c}^{-1}\Theta_{c}(z)+I
\end{array}
\right),
\end{align*}
} here  \begin{align*} \Delta_{l,l}=&\mbox{\emph{diag}}\{\sigma_r,
r=\overline{1,
R}\}\oplus \mbox{\emph{diag}}\{\lambda^{(1)}_w, w=\overline{1, W}\} \\
&\oplus \mbox{\emph{diag}}\{\lambda^{(2)}_v, v=\overline{1,
V}\}\oplus [\mbox{\emph{diag}}\{s_m, m=\overline{1, M}\}]^{\oplus
l},  \ \ \ l=\overline{g, c},
\end{align*}
\begin{align*}
\Gamma_{l,l'}= I_{RW}\otimes E_r\otimes
I_{M^l}\otimes\varsigma^{\otimes r},  \ \ \ l'=l+r, r=\overline{1,
c-1-l}, l=\overline{g, c-2}.
\end{align*}

It is seen from \eqref{38} that the matrix $Y(z)$ is reducible and
the matrix $Y(1)$ contains only one irreducible stochastic diagonal
block $Y_{22}(1)$. According to the Proposition 1.2, we get a
sufficient condition for ergodicity of the Markov chain $\{\xi_t, t
\geq0\}$ which is the fulfillment of the following inequality:
\begin{align}\label{39}
\left[\det(zI-Y_{22}(z))\right]'|_{z=1}>0.
\end{align}

A more constructive form of ergodicity condition \eqref{39} in our
system is given by the following theorem.
\begin{theorem}
A sufficient condition for ergodicity of the Markov chain $\{\xi_t,
t \geq0\}$, is the fulfillment of the inequality
\begin{align*}
\frac{\lambda_1}{\bar{\mu}_1}+\frac{\lambda_2}{\bar{\mu}_2}<1,
\end{align*}
where
\begin{align*}
\bar{\mu}_1=\mathbf{X}_1S_0^{\oplus g}\mathbf{e}_{M^{g-1}}, \ \ \
\bar{\mu}_2=\mathbf{X}_2S_0^{\oplus c}\mathbf{e}_{M^{c-1}},
\end{align*}
and $\mathbf{X}_i$, $i=1,2$, are the unique solutions of the
following system of linear algebraic equations:
$$
\left\{
\begin{array}{l}
 \mathbf{X}_1[S^{\oplus g}+S_0^{\oplus
g}(I_{M^{g-1}}\otimes
\varsigma)]=0,\\
\mathbf{X}_1\mathbf{e}=1,
\end{array}
\right.
$$
$$\left\{
\begin{array}{l}
\mathbf{X}_2[S^{\oplus c}+S_0^{\oplus c}(I_{M^{c-1}}\otimes
\varsigma)]=0,\\
\mathbf{X}_2\mathbf{e}=1.
\end{array}
\right.
$$
\end{theorem}
\begin{proof}
The proof of this theorem follows the steps given in the paper
 \cite{Breuer2002} and consequently we omit the details
here.
\end{proof}

 Next, we assume that the stability condition is fulfilled.
Denote the steady state probabilities of the Markov chain as
\begin{align*}
&P\left(o, b, \nu, \gamma, m^{(1)},\cdots, m^{(b)}\right)\\
&=\lim_{t\to\infty} P\left\{o_t=o, b_t=b, \nu_t=\nu,
\gamma_t=\gamma,
m_t^{(1)}=m^{(1)},\cdots,m_t^{(b_t)}=m^{(b)}\right\}.
\end{align*}
Let the row-vector $\mathbf{P}_i$ of dimension $K$ denote the
stationary-state probabilities $$\mathbf{P}_i=P\left(i, b, \nu,
\gamma, m^{(1)},\cdots, m^{(b)}\right), \ \ \ i\geq0,$$ and define
the infinite-dimensional probability vector
$\mathbf{P}=(\mathbf{P}_1, \mathbf{P}_2, \cdots)$. When the system
is stable, $\mathbf{P}$ is the unique solution to $\mathbf{P}
\mathbf{Q}=0$ and $\mathbf{P} \mathbf{e}=1$. However, it is still an
open problem to express $\mathbf{P}$ in the closed form (e.g.,
generating function), and even though $\mathbf{Q}$ is highly
structured, $\mathbf{P}$ cannot be expressed in a tractable
analytical form. In the following, we will adopt the algorithm of
asymptotically quasi-Toeplitz Markov chains to solve the equilibrium
equation $\mathbf{P} \mathbf{Q}=0$.

\section{Algorithm}\label{sec4}
In this section, we use the numerically stable algorithm which has
been elaborated in \cite{Klimenok2006} to calculate the stationary
distribution of the system. The algorithm is based on censoring
technique of the asymptotically quasi-Toeplitz Markov chain and
consists of some essential steps which are given in the following
Algorithm 4.1.

\begin{algorithm}
Computing the probability vectors $\mathbf{P}_i$:\\
Step 1:  Calculate the matrix $G$ as the minimal nonnegative
solution to the non-linear matrix equation:
\begin{align}\label{eq:41}
G=Y(G).
\end{align}
We should note that the first $RWV\frac{1-M^g}{1-M}$ rows of the
matrix $G$ coincide with corresponding rows of the the matrix $Y(1)$
and the rest $RWV\sum_{n=g}^cM^n$ entries of the matrix $G$ can be
calculated by the iterative method which is available in \cite{Lucantoni}.\\
Step 2: Chose the integer $k_0$ as a minimal value of $k$ for which
the norm $\|G-\sum_{n=k}^\infty Q_{k+1,n}G^{n-k}\|$ of the residual
is less than the preassigned value $\epsilon$.\\
Step 3:  Calculate the matrices $G_{k_0-1}$, $G_{k_0-2}$ , $\cdots$,
$G_{0}$ from the following recursive equation
\begin{align}
G_k=\left(-\sum_{j=k+1}^\infty Q_{k+1,n}G_{j-1}G_{j-2}\cdots
G_{k+1}\right)^{-1}Q_{k+1,k},
\end{align}
$k=0, 1, \cdots, k_0-1$, with the boundary condition $G_k=G$, $k\geq
k_0$.\\
Step 4: Calculate the matrices $\bar{H}_{i,j}$, $j\geq i$,
$i\geq0$,using the formulae
\begin{align}
\bar{H}_{i,j}=H_{i,j}+\sum_{k=j+1}^\infty
H_{i,k}G_{k-1}G_{k-2}\cdots G_j, \ \ \ j\geq i, \ \ \  i\geq0,
\end{align}
where $G_k=G$, $k\geq
k_0$.\\
 Step 5: Calculate the matrices $F_j$, $j\geq0$, according to
the recursive relation
\begin{align}
F_0=I, \ \ \
F_j=\sum_{i=0}^{j-1}F_i\bar{H}_{i,j}(-\bar{H}_{j,j})^{-1}, \ \ \
j\geq1.
\end{align}
Step 6: Calculate the vector $\mathbf{P}_0$ as the unique solution
to the system of the linear equations
\begin{align*}
\left\{
\begin{array}{l}
\mathbf{P}_0(-\bar{H}_{0,0})=\mathbf{0},\\
\mathbf{P}_0\sum\limits_{j=0}^\infty F_j\mathbf{e}=1.
\end{array}
\right.
\end{align*}
Step 7: Calculate the vectors $\mathbf{P}_j$, $j\geq1$, as
\begin{align*}
\mathbf{P}_j=\mathbf{P}_0F_j, \ \ \ j\geq1.
\end{align*}
\end{algorithm}

\begin{remark}
We can eliminate the infinite sums in Step 6 as the following
method: pick $N>k_0$, and replace the second equation in Step 6 with
$\mathbf{P}_0\sum\limits_{j=0}^N F_j\mathbf{e}=1$, such that the
inequality $\mathbf{P}_N\mathbf{e}<\epsilon_0$ is satisfied, where
$\epsilon_0$ is some preassigned small number. Note that the inverse
matrices that appear in the algorithm exist and are nonnegative. So
the computation of the algorithm is stable. The algorithm can be
realized in framework of extension of ``SIRIUS++" software
\cite{SIRIUS} that was developed using the results of the
asymptotically quasi-Toeplitz Markov chains. As such the
determination of the stationary distribution is straightforward.
\end{remark}

\section{performance measures}
Having the stationary distribution $\mathbf{P}_i, i\geq0$ been
calculated, we are able to calculate  some stationary performance
measures of
the system under consideration in this section.\\
\indent Let $\mathbf{P}_i=(\mathbf{P}_{i0}, \mathbf{P}_{i1}, \cdots,
\mathbf{P}_{ic})$, $i\geq0$, where $\mathbf{P}_{ib}$ of size
$1\times RWVM^b$ is a row vector, $0\leq b\leq c$. Denote $P(i,b)$
the joint probability to have $i$ repeated customers in the orbit
and $b$ customers on the servers at arbitrary time. From the results
of the above section, we get
\begin{align*}
P(i,b)=\mathbf{P}_{ib}\mathbf{e}&=\sum_{r=1}^R\sum_{\nu=1}^W\sum_{\gamma=1}^V
\sum_{m^{(1)}=1}^M\cdots\sum_{m^{(b)}=1}^MP\left(i,
b, r, \nu, \gamma, m^{(1)},\cdots, m^{(b)}\right)\\
&=\mathbf{P}_i\mathbf{e}\left[WV\frac{1-M^b}{1-M}+1,
WV\frac{1-M^{b+1}}{1-M}\right],
\end{align*}
where $\mathbf{e}\left[RWV\frac{1-M^b}{1-M}+1,
RWV\frac{1-M^{b+1}}{1-M}\right]$ of size $K\times1$ is a column
vector of suitable size having ones as its
$RWV\frac{1-M^b}{1-M}+1$th to $RWV\frac{1-M^{b+1}}{1-M}$th entries
and zeros as the rest entries.

\begin{corollary}
 \ \
\begin{enumerate}
\item The probability that there are $i$ repeated customers in the orbit
at arbitrary time
$$P(i,\bullet)=\sum\limits_{b=0}^{c}P(i,b)=\mathbf{P}_i\mathbf{e},
 \ \ \ i\geq0.$$
\item  The probability that there are $b$ busy servers at
arbitrary time
$$P(\bullet,b)=\sum\limits_{i=0}^\infty P(i,b),
 \ \ \ b=\overline{0, c}.$$
\item  The mean number of busy servers
$$L_b=\sum\limits_{b=1}^{c}bP(\bullet,b).$$
\item  The mean number of repeated customers in the orbit
$$L_{orb}=\sum\limits_{i=1}^\infty i\mathbf{P}_i\mathbf{e}.$$
\item  The mean number of customers (primary customers and priority customers) in
the system
$$L_s=L_{orb}+L_b.$$
\item The blocking probability for an arbitrary primary customer
$$P_{b1}=1-\frac{1}{\lambda_1}\sum_{n=1}^g\sum_{i=0}^\infty P(i,g-n)\sum_{k=0}^n(k-n)(I_R\otimes D_k\otimes
I_{VM^{g-n}})\mathbf{e}.$$
\item The blocking probability for an arbitrary batch of primary customers
$$P_{bb1}=1-\frac{1}{\lambda_{b1}}\sum_{n=1}^g\sum_{i=0}^\infty P(i,g-n)\sum_{k=1}^n(I_R\otimes D_k\otimes
I_{VM^{g-n}})\mathbf{e}.$$
\item The blocking probability for an arbitrary priority customer
$$P_{b2}=1-\frac{1}{\lambda_2}\sum_{n=1}^c\sum_{i=0}^\infty P(i,c-n)\sum_{k=0}^n(k-n)(I_{RW}\otimes E_k\otimes
I_{M^{c-n}})\mathbf{e}.$$
\item The blocking probability for an arbitrary batch of priority customer
$$P_{bb2}=1-\frac{1}{\lambda_{b2}}\sum_{n=1}^c\sum_{i=0}^\infty P(i,c-n)\sum_{k=1}^n(I_{RW}\otimes E_k\otimes
I_{M^{c-n}})\mathbf{e}.$$
\item Busy period: If we define the busy period $B$ of this queueing system with
repeated customers as the period that starts at the epoch when an
arriving batch of customers (primary or priority) finds an empty
system (all servers are idle and no repeated customer in the orbit)
and ends at the departure epoch at which the system is empty again,
then the mean busy period is given by
$$E(B)=\frac{1}{\lambda_1+\lambda_2}\left(\frac{1}{P(0,0)}-1\right).$$
\end{enumerate}
\end{corollary}

Moreover, using the approach based on Markov renewal processes, we
can obtain the following statement.
\begin{proposition}
The generating function $\Pi(z)$ of the queue state stationary
distribution at the departure epochs is characterized by the
following functional-differential equation:
\begin{align*}
\Pi'(z)=\Pi(z)A(z),
\end{align*}
where
\begin{align*}
A(z)=B^{-1}(z)B'(z)-B^{-1}(z)\Delta (T_1^{-1}\otimes
I_{WV\sum_{k=0}^cM^k})z^{-1}B(z)+z^{-1}(T_1^{-1}\otimes
I_{WV\sum_{k=0}^cM^k})B(z),
\end{align*}
\begin{align*}
B(z)=[\Upsilon-z\hat{I}-\Delta^{-1}\hat{I}z\Gamma(z)][\Upsilon-z\hat{I}-\Delta^{-1}\bar{I}z\Gamma(z)]^{-1}.
\end{align*}
\end{proposition}

\begin{proof}
The proof can be implemented similarly as the proof of
\cite{Breuer2002} and is omitted here.
\end{proof}

\section{Optimization}
In this section, we borrow some ideas and methods from
\cite{Trivedi} to investigate the optimization problem that how to
choose the optimal values of $g$ and $c$ so as to minimize the
priority customer blocking probability as well as the primary
customer blocking probability. Therefore, we obtain a
multi-objective optimization problem in which the number of total
servers $c$ and the parameter $g$ which is equal to the number of
non-guard servers are the decision variables. However, due to the
complexity, we may fix $c$ and only consider $g$ as the decision
variable to simply the optimization problem. In the case of the two
objectives, we can choose either $P_{b1}$ or $P_{b2}$ as the
objective function to be minimized and impose a constraint on the
other one. If we regard the two blocking probabilities as two binary
functions of $g$ and $c$ two variables, then we write two
representative optimization problems as follows:
\begin{enumerate}[(I)]
\item  Given $T_k, D_k, E_k, k=0, 1, \cdots, S$ and $c$ determine the
optimal integer values of $g$ such that
\begin{align*}
&\mbox{Mininize} \ \ \ P_{b1}(g)\\
& \mbox{subject to} \ \ \  P_{b2}(g)\leq p_0.
\end{align*}

\item Given $T_k, D_k, E_k, k=0, 1, \cdots, S$  determine the
optimal integer values of $g$ and $c$ such that
\begin{align*}
&\mbox{Mininize} \ \ \ c\\
& \mbox{subject to} \ \ \ P_{b1}(g,c)\leq p_1, \ \ \ \mbox{and} \ \
\ P_{b2}(g,c)\leq p_2.
\end{align*}
\end{enumerate}
Here, the constants $p_0$,  $p_1$ and  $p_2$ are three preassigned
numbers.

We first consider the optimization problem (I). It is clear that $
P_{b2}(g)> P_{b2}(g-1)$ (see Fig. 3). Thus, we can get the largest
value of $g$ such that $P_{b2}(g)\leq p_0$. Then based on the
monotonicity property $P_{b1}(g)> P_{b1}(g+1)$, we know that such a
value of $g$ will minimize $P_{b1}(g)$. Therefore, we can obtain the
optimal value of $g$ using a simple one dimensional search over the
range $\{1,\cdots, c-1\}$ for $g$ such that
$g^*=\max\{g|P_{b2}(g)\leq p_0\}$.

Then we study the optimization problem (II). Since we cannot derive
the analytical solution of this system, in the next section, we will
resort to the numerical examples to illustrate the optimization
problem (II). Here, we provide an algorithm to obtain the minimum
value of $c$.
\begin{algorithm} Solving the optimization problem (II)\\
Step 1: Set the values $p_1$ and  $p_2$ and let c=2;\\
Step 2: For $g=1$ to $c-1$, find the feasible region for $g$:\\
$$\Theta=\{g|P_{b1}(g,c)\leq p_1, \ \ \ \mbox{and} \ \ \
P_{b2}(g,c)\leq p_2\};$$ Step 3: If $\Theta\neq\emptyset$, then stop
and obtain $c^*=c$. Otherwise,
put $c=c+1$ and goto step 2; \\
Step 4: Obtain $c^*=c$.
\end{algorithm}

\section{Application to cellular wireless networks and numerical examples}
In this section, we present an application to cellular wireless
networks for the model under discussion and some numerical results
that show the effect of the system parameters on the performance
measures and channel allocation schemes in the cellular mobile
networks. In the cellular mobile wireless network, the regional
service area is divided into multiple adjacent cells, each of them
served by a base station (BS) with a limited number of channels.
Here, we assume the number of channels is equal to $c$, $c\geq2$. A
mobile subscriber who has wireless terminal such as mobile phone for
voice service, wireless data terminal for data service, dual
terminal for voice and data, and so on, communicates via radio links
to a BS, one for each cell. Generally, there are two types of calls
in a specific cell in the cellular mobile networks: originating
calls, i.e., calls originating in that cell, and handoff calls,
i.e., current ongoing calls caused by mobile users moving from one
cell to that cell. We assume that the originating calls arrive at
the system according to a BMAP with the characterizing matrix
sequence $(D_n: n\in \mathbb{N}_0)$ and the handoff calls arrive at
the system according to a BMAP with the characterizing matrix
sequence $(E_n: n\in \mathbb{N}_0)$.  When a batch of new calls
originates in that cell, some idle channels are required to the BS
for the new calls setup. If there are idle channels in the BS of
that cell, then the BS assigns some idle channels to the originating
calls; otherwise some new calls (maybe all the batch) in that batch
are blocked and the blocked mobile subscribers may try their luck
later. Similarly, if some idle channels are available in the BS of
that cell, the batch of handoff calls is successfully handed over
without an interruption; otherwise some handoff calls (or all the
batch) are blocked and retry for service after a short period. Any
blocked call (originating or handoff) retries for service where the
rates of the retrials depend on a changing environment. Thus, we use
the MMPP with the characterizing matrices $(T_0, T_1)$ to model the
retrial process. The quality of service is disturbed seriously if a
handoff call becomes blocked when it crosses the cell boundary.
Generally, the base station gives priority to a handoff call over an
originating call. The usual way to do this is the scheme with guard
channels, see \cite{Guerin} and
 \cite{Hong}. Here, we assume the number of
non-guard channels is equal to $g$. Thus, the number of guard
channels is equal to $c-g$. Furthermore, we suppose the conversation
time of an ongoing call (originating or handoff) by a channel is of
PH type with a representation $(\varsigma, S)$. Then the cellular
mobile wireless network can be modeled as a $BMAP_1, BMAP_2/PH/g, c$
retrial queueing system where the rate of individual repeated
attempts from the orbit is modulated according to a MMPP. The
originating calls can be regarded as primary customers and the
handoff calls can be regarded as priority customers.

In the following, we give some simple numerical examples that
illustrate the effects of various parameters on the system
performance measures of the cellular mobile wireless network. The
algorithms 4.1 and 6.1 developed in sections 4 and 6 have been
written into some MATLAB programs.  Realization of these algorithms
on computer does not meet any difficulty. Only the calculation time
can be long especially when the dimension of the system is large.
For example, if we let $c=8$, then the dimension $K=4088$, and the
calculation time is 98s. For $c=100$, the dimension is equal to
1.0141e+031, and the calculation time is very long. Thus, the
calculation time increases drastically with the increase of the
number of servers $c$. Note that in all the below examples, we
choose the parametric values in a way such that the system is
stable. Numerical results are showed in Tabs. 1-3 and Figs. 1-3. For
comparison, in each of the pictures, we also plot three curves which
correspond to $\sigma=6.6429$, $13.2857$ and $26.5714$.

For the convenience of the numerical calculation, here, we suppose
that handoff calls and originating calls arrive to the system
according to two independent special BMAPs (In fact, they are MAPs)
which are characterized by the following matrices:
$$D_0=\lambda_o\left(\begin{array}{cc}
-11&2\\
5&-20
\end{array}
\right),  \ D_1=\lambda_o\left(\begin{array}{cc}
8&1\\
3&12
\end{array}
\right),  \ D_k=O_2, \ k\geq2,$$
$$E_0=\lambda_h\left(\begin{array}{cc}
-3&0\\
1&-2
\end{array}
\right),  \ E_1=\lambda_h\left(\begin{array}{cc}
1&2\\
0&1
\end{array}
\right),  \ E_k=O_2, \ k\geq2,$$
 This means the stationary arrival
rates of the originating calls and the handoff calls are
$\lambda_1=10.6364\lambda_o$ and $\lambda_2=1.6667\lambda_h$,
respectively.

The repeated calls in the orbit repeat there attempts to reach a
server according to the MMPP which is defined by the matrices:
$$T_0=\lambda_r\left(\begin{array}{cc}
-15&3\\
4&-19
\end{array}
\right), \ \ \ T_1=\lambda_r\left(\begin{array}{cc}
12&0\\
0&15
\end{array}
\right), $$ Thus, the mean retrial rate for the mobile engaged in
the cell is $\sigma=13.2857\lambda_r$.

 The PH service process is defined by the matrices:
$$S=\left(\begin{array}{cc}
-23&9\\
14&-17
\end{array}
\right), \ \ \ \varsigma=(0.4, 0.6).$$ So that the mean rate of
service is $\mu=8.1288$.

\indent From the above descriptions, we have $ W=2, V=2, M=2$ and
$R=2$.
\begin{center}
\begin{table}
\footnotesize \caption{The stationary join distribution of the
system with $(c, g, \lambda_o, \lambda_h,  \lambda_r)=(8, 6, 2, 2,
2)$}
    \begin{tabular}{lcccccccccc}

        \hline
       $i$ $\backslash$ $b$ & 0&1 & 2&3 & 4 & 5 & 6& 7&8&$sum$\\
        \hline
        0& 0.0467  &  0.1413  &  0.2136  &  \textbf{0.2148}  &  0.1604   & 0.0923   & 0.0376  &  0.0016 &   0.0001&  \textbf{0.9084}\\
        1& 0.0001  &  0.0006  &  0.0020  &  0.0047  &  0.0091   & 0.0149   & 0.0208   & 0.0013  &  0.0001 & 0.0536\\
        2& 0.0000  &  0.0000  &  0.0002  &  0.0008  &  0.0023   & 0.0054   & 0.0109   & 0.0008   & 0.0000 & 0.0206\\
        3&  0.0000 &   0.0000 &   0.0000 &   0.0002 &   0.0007  &  0.0022  &  0.0055  &  0.0005  & 0.0000 & 0.0092\\
        4& 0.0000  &  0.0000  &  0.0000  &  0.0001  &  0.0002   & 0.0009   & 0.0028   & 0.0003   & 0.0000 & 0.0043\\
        5& 0.0000  &  0.0000  &  0.0000  &  0.0000  &  0.0001   & 0.0004   & 0.0014   & 0.0001   & 0.0000& 0.0020\\
        & 1.0e-003 $\times$\\
        6& 0.0000  &  0.0001  &  0.0007  &  0.0055  &  0.0334   & 0.1661   & 0.6986   & 0.0667   & 0.0051 & 0.9761\\
        &1.0e-003 $\times$\\
        7& 0.0000  &  0.0000   & 0.0002  &  0.0019  &  0.0130   & 0.0731   & 0.3474   & 0.0337   & 0.0027 & 0.4721\\
        & 1.0e-003$\times$\\
        8&0.0000  &  0.0000   & 0.0001  &  0.0007   & 0.0052   & 0.0325   & 0.1724   & 0.0169   &
        0.0014 & 0.2291\\
        & 1.0e-004$\times$\\
        9&0.0000  &  0.0000   & 0.0002  &  0.0025 &   0.0211  &  0.1460   & 0.8533   & 0.0846  &
        0.0069& 1.1145\\
        &1.0e-004$\times$\\
        10 &0.0000  &  0.0000   & 0.0001  &  0.0009  &  0.0087 &   0.0660  &  0.4217  &  0.0421  &
        0.0035& 0.5429\\
\hline
 $sum$ &  0.0468 &  0.1419 & 0.2158 & \textbf{0.2206}  & 0.1728 & 0.1162&  0.0800 & 0.0047& 0.0002& 0.999\\
        \hline
    \end{tabular}
    \end{table}
\end{center}
\begin{table}
 \caption{The the optimal value $g^*$ for the
optimization problem (I) with $(c, \lambda_o, p_0)=(20, 10,
0.0001)$}
    \begin{tabular}{lcccccccccccccccc}

        \hline
        $\lambda_r$ $\backslash$ $\lambda_h$ & 1&2 & 3&4 & 5 & 6 & 7& 8&9&10& 11 & 12& 13&14&15&20\\
        \hline
        1& 18 & 18 & 18 & 17  &  17  &  17 &  16& 16&16&16&  15 &  15& 15&14&14&13\\
        10&18 & 18 & 18 & 17  &  17  &  17 &  16& 16&16&16&  15 &  15& 15&14&14&13\\
        20& 18 & 18 & 18 & 17  &  17  &  17 &  16& 16&16&16&  15 &  15& 15&14&14&13\\
        \hline
    \end{tabular}
    \end{table}

\begin{table}
 \caption{The optimal value $c^*$ for the
optimization problem (II) with $(\lambda_r, p_1, p_2)=(10, 0.001,
0.0001)$}
    \begin{tabular}{lcccccccccccccccc}

        \hline
        $\lambda_h$ $\backslash$ $\lambda_o$ & 1&2 & 3&4 & 5 & 6 & 7& 8&9&10& 11 & 12& 13&14&15&20\\
        \hline
        1& 8 & 11 & 14 & 16  &  18  &  20 &  22& 24&26&28&  29 &  31& 33&35&37&45\\
        5&10 & 13 & 15 & 17  &  19  &  21 &  23& 25&27&29&  30 &  32& 34&36&38&46\\
       10& 13 & 15 & 17 & 19  &  21  &  23 &  25& 27&29&31&  32 &  34& 36&37&39&47\\
        \hline
    \end{tabular}
    \end{table}

\begin{figure}
\centering{}\includegraphics[width=4.5in]{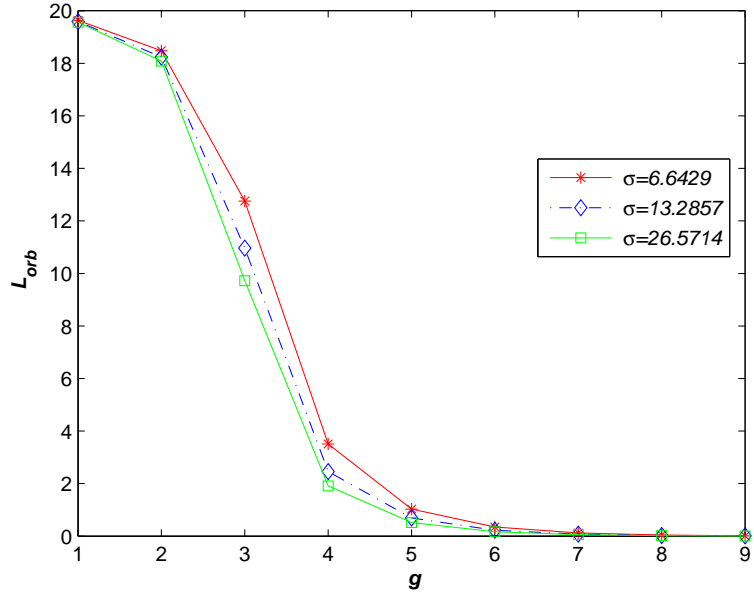}\caption{$L_{orb}$
as function of the parameter $g$ with $(c, \lambda_o,
\lambda_h)=(10, 2, 2)$.}
\end{figure}
\begin{figure}
\centering{}\includegraphics[width=4.5in]{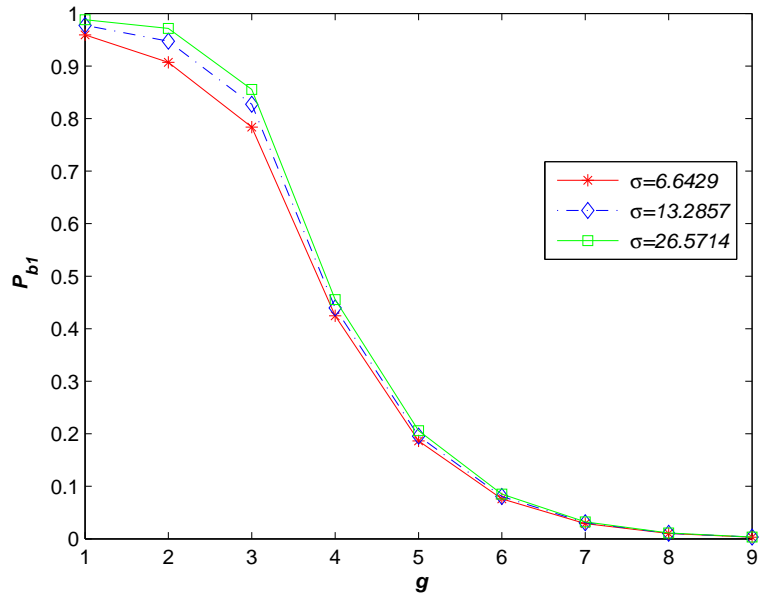}\caption{Dependence
of the blocking probability for originating calls on the value $g$
with $(c,\lambda_o, \lambda_h)=(10,2,2)$.}
\end{figure}
\begin{figure}
\centering{}\includegraphics[width=4.5in]{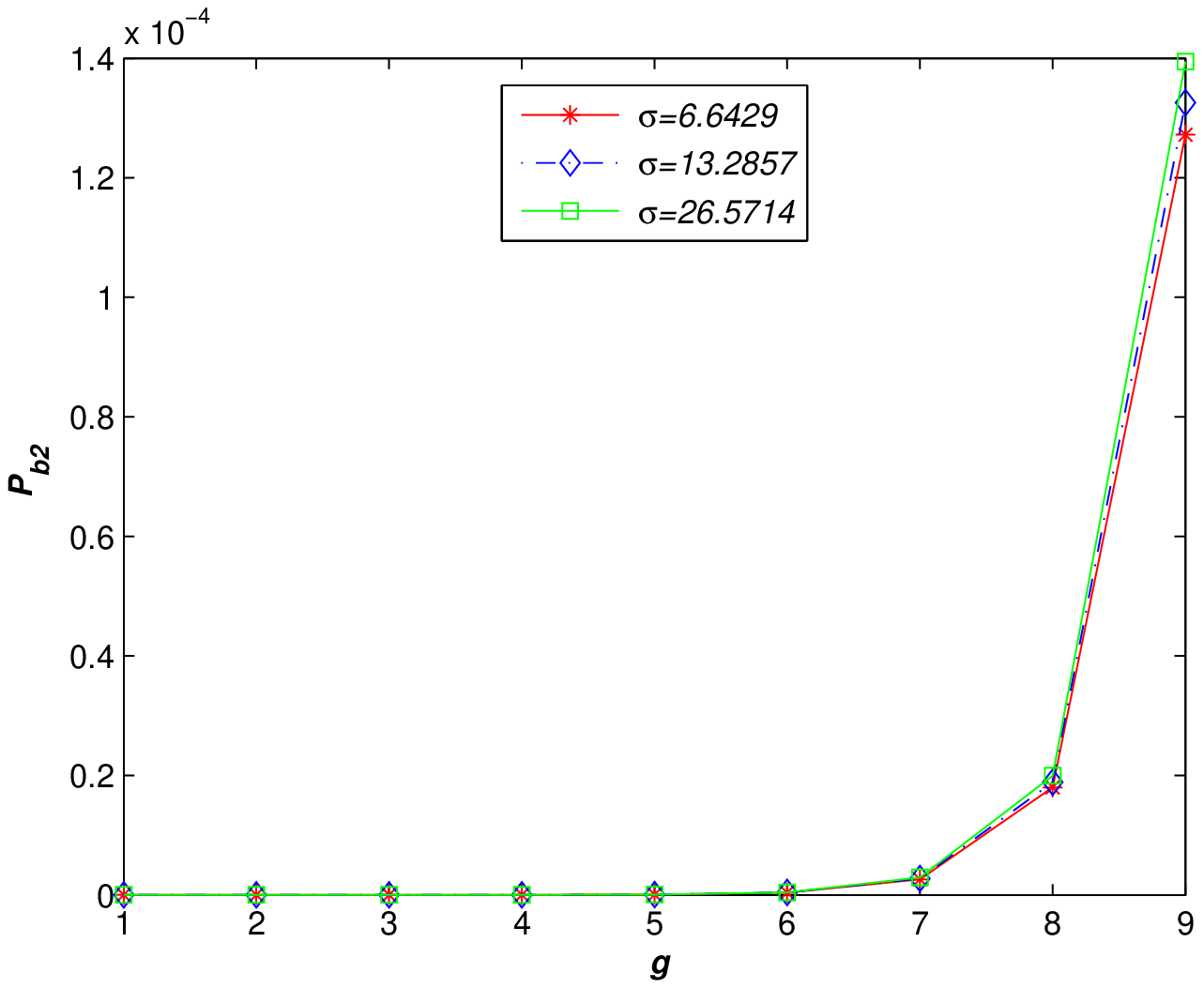}\caption{Dependence
of the blocking probability for handoff calls on the value $g$ with
$(c,\lambda_o, \lambda_h)=(10,2,2)$.}
\end{figure}

The joint distribution of the number of busy channels and the number
of repeated calls in the orbit is presented in Tab. 1 with $(c, g,
\lambda_o, \lambda_h)=(8, 6, 2, 2)$. From Tab. 1, we observe that
the joint probability that there are three busy channels and zero
repeated call  at arbitrary time is equal to 0.2148 which is the
maximum value and the probability that there is no repeated call in
the orbit is equal to 0.9084 which is the maximum value of the
marginal probabilities (lateral). The probability that there are
three busy channels  at arbitrary time  is equal to 0.2206 which is
the maximum value of the marginal probabilities (longitudinal).

The effects of the retrial rate and the arrival rate of the
originating calls on the optimal value of $g$ for the optimization
problem (I) with the set of parameters $(c, \lambda_h,
\theta_0)=(20, 10, 0.0001)$ are reported in Tab. 2. We observe that
the optimal value $g^*$ decreases monotonously as the arrival rate
of the handoff calls increases. From this, we can predict that the
optimal value $g^*$ increases monotonously as the arrival rate of
the originating calls increases. However, the retrial rate of the
repeated calls in the orbit has no significant impact on the optimal
value $g^*$. The appearance of this phenomenon may be based on the
fact that the mean number of repeated calls is very small compared
to the mean number of arriving calls. Hence it may be concluded that
the arrival rates of the two types of calls are the major factors in
deciding the optimal numbers of guard channels.

Tab. 3 shows the effects of the arrival rates of the originating
calls and the handoff calls on the optimal value $c^*$ for the
optimization problem (II), where we set $(\lambda_r, p_1, p_2)=(10,
0.001, 0.0001)$. As is to be expected, when the arrival rates of the
two types of calls increase, the optimal value $c^*$ increases.
Moreover, we observe that the optimal value $c^*$ increases faster
with the arrival rates of the originating calls than that with the
arrival rates of the handoff calls. Therefore, we understand that
the arrival of originating calls plays a most significant role in
deciding the optimal value $c^*$. In other word, this type of calls
whose arrival rate is largest will most strongly affect the optimal
design of the system.

Fig. 1 illustrates the dependence of the mean number of repeated
calls in the orbit on the value of parameter $g$ where we set
$(c,\lambda_o,\lambda_h)=(10, 2, 2)$. As it can be seen from Fig. 1
that $L_{orb}$ decreases monotonously as the value $g$ increases,
especially when $g\leq6$. But when $g\geq7$, the curves of the
values $L_{orb}$ are almost overlapping with $x$-axis. That is to
say, there are very few repeated customers in the orbit when the
number of guard channels is small.

Figs. 2 and 3 show the dependence of the blocking probability for an
arbitrary originating call and handoff call on the value of
parameter $g$ respectively when $(c,\lambda_o,\lambda_h)=(10, 2,
2)$.  From Figs. 2 and 3, we find that the blocking probability for
an arbitrary originating call $P_{b1}$ decreases monotonously, but
the blocking probability for an arbitrary handoff call $P_{b2}$
increases monotonously as the value $g$ increases, which agrees with
the intuitive expectations. In addition, when $g\leq6$, the curves
of the values $P_{b2}$ are almost overlapping with $x$-axis. It
means that the blocking probability for an arbitrary handoff call
$P_{b2}$ is roughly equal to zero when the number of guard channels
is large.

We must note that, in each of the pictures, the three curves are
almost overlapping. This phenomenon tells us that the retrial rate
$\sigma$ has little impact on the system performance measures.

From what has been described above, we can draw the conclusion that
the performance measures are mainly affected by the two types of
calls' arrivals and service patterns, but the retrial rate plays a
less crucial role.

 \section{Conclusions}\label{sec8}
 In this paper, we have investigated the BMAP/PH/N type retrial queue
 with two types of customers. Our work can be considered as an extension of
  \cite{Breuer2002} and \cite{Dudin2012}. The behavior of this model is
 described by a multi-dimensional continuous-time asymptotically
quasi-Toeplitz Markov chain. Sufficient condition for the ergodicity
of the Markov chain is given. An algorithm for computing the
stationary distribution is presented. Expressions for calculation of
main performance measures of the model are derived. The optimization
problem how to choose the optimal values of $g$ and $c$ is discussed
and an algorithm is also provided. An application of the model to
the cellular wireless network is implemented. The dependences of the
main performance characteristics of the model on its parameters are
graphically demonstrated. We obtain the results that the retrial
pattern has very little impact on the optimal values of guard
channels and total channels, the number of busy channels and the
blocking probabilities for the two types of calls.

Future work will consider the BMAP/PH/c retrial queues with
vacations or working vacations and their applications to computer
communication networks.

\section*{Compliance with Ethical Standards}
 This study was funded by the National Natural Science
Foundation of China (11201489, 11271373, 11371374). The authors
declare that they have no conflict of interest.


\bibliographystyle{model5-names}
\bibliography{BMAP}

\end{document}